\documentclass[12pt]{article}

\usepackage{amssymb,amsmath,latexsym}
\usepackage[dvips]{graphics}

\input epsf.tex

\newtheorem{thm}{Theorem}[section]

\newtheorem{remark}[thm]{Remark}
\newtheorem{defn}[thm]{Definition}
\newtheorem{corollary}[thm]{Corollary}

\newtheorem{lemma}[thm]{Lemma}

\numberwithin{equation}{section}

\def\RR{{\bf R}}
\def\ZZ{{\bf Z}}
\def\DD{{\bf D}}
\def\EE{{\cal E}}

\def\RR{{\bf R}}

\def\wt{\widetilde}
\def\pf{\noindent{\bf Proof.} }
\def\eps{{\varepsilon}}
\def\E{{\bf E}}

\def\P{{\bf P}}

\def\dist{{\hbox{\rm dist}}}
\def\qed{{\hfill $\Box$ \bigskip}}

\begin{document}
\title{\bf Boundary Trace of Reflecting Brownian Motions
}
\author{{\bf Itai Benjamini}, \quad  {\bf Zhen-Qing Chen}\thanks{This research
is supported in part by NSF Grant
DMS-0071486 and a RRF grant from University of Washington.} \quad
and \quad {\bf Steffen Rohde}\thanks{This research
is supported in part by NSF Grant
DMS-0201435.}
}
\bigskip
\date{(July 14, 2003) }
\maketitle

\bigskip

\begin{abstract}
We establish a uniform dimensional result for normally reflected
Brownian motion (RBM) in a large class of non-smooth domains.
Exact Hausdorff dimensions for the boundary occupation time
and the boundary trace of RBM are given.
Extensions to stable-like jump processes and to
symmetric reflecting diffusions are also mentioned.
\end{abstract}

\vspace{.6truein}


\bigskip


\vfill \eject

\begin{doublespace}

\section{Introduction}

Let $n\geq 2$ and
$D\subset \RR^n$ be a domain (connected open set) with compact closure.
Consider a reflecting Brownian motion (RBM in abbreviation)
$X$ in $D$. Heuristically, RBM in $D$ is a continuous Markov process
$X$ taking values in $\overline D$ that behaves like a Brownian motion
in $\RR^n$ when $X_t\in D$ and is instantaneously pushed back
along the inward normal direction when $X_t\in \partial D$.
RBM on smooth domains can be constructed in various ways,
including the deterministic Skorokhod problem method, stochastic
differential equation with boundary condition, martingale problem
methods, etc. see the Introduction of \cite{C}.
When $D$ is non-smooth, all the above mentioned methods cease to
work.
On non-smooth domains, the most effective way to construct
RBM is to use the Dirichlet form method, which will be recalled
in Section \ref{S:2}. The RBM constructed through
Dirichlet form coincides with all the other standard definitions
in smooth domains. Using the Dirichlet form approach, it can be shown
that, when $D$ is a simply connected planar domain, RBM $X$ in $D$ is
the time change of
the conformal image of RBM in a unit disc.

\medskip

This paper is concerned with  Hausdorff dimensions
of various random sets associated with RBM.
The study of Hausdorff dimensions of random sets associated
with Brownian motion, stable processes and more generally
L\'evy processes together with their fractal structures
has been an active research area in the last 40 years.
See Xiao \cite{Xiao} for a recent survey on this subject.
However, to the authors' knowledge,
this is the first time that such a study has been
conducted for RBM in Euclidean domains.

\medskip

The results obtained in this paper hold for a large class
of non-smooth domains. In order to convey our results
as transparently as possible, in this introductory section
we confine ourselves to a special case
of the more general results established in this paper by
assuming
that $D$ is a bounded {\it uniform} domain.
The following definition is taken from V\"ais\"al\"a \cite{Va},
where various equivalent definitions are discussed.

\begin{defn}\label{D:2.1}
A domain $D\subset \RR^n$ is called uniform if there exists a constant $C$
such that for every $x,y\in D$ there is a rectifiable curve $\gamma$
joining $x$ and $y$ in $D$ with
$\text{length}(\gamma)\leq C|x-y|$ and moreover
$\min \left\{ |x-z|, \, |z-y|)
\right\}
\leq C \, \dist(z,\, \partial D)$
for all points $z\in\gamma$.
Here $\dist(z,\, \partial D)$ is the Euclidean distance
between point $z$ and the set $\partial D$.
\end{defn}

For example, the classical van Koch snowflake
domain in the conformal mapping theory is a uniform domain
in $\RR^2$.
The uniform domain is also called $(\eps, \infty)$-domain
in the terminology of Jones \cite{J}.
Note that every Lipschitz domain is uniform, and every
{\it non-tangentially accessible domain} defined by Jerison and
Kenig in \cite{JK} is a  uniform domain (see (3.4) of \cite{JK}).
However, the boundary of a uniform domain can be highly
nonrectifiable and, in general, no regularity of its boundary
can be inferred (besides the easy fact that the Hausdorff dimension
of the boundary is strictly less than $n$).
For any $\alpha \in [n-1, n)$,
one can construct a uniform domain
$D\subset \RR^n$ such that ${\cal H}^\alpha (U \cap \partial D)>0$
for any open set $U$ satisfying $U\cap \partial D \not= \emptyset$.
Here ${\cal H}^\alpha$  denotes the $\alpha$-dimensional
Hausdorff measure in $\RR^n$.

\medskip

In this paper,
we first extend Kaufman's uniform doubling dimension result
\cite{Ka} for planar Brownian motion to
RBMs in bounded domains. It follows as a special case
of Theorem \ref{T:2.1} below
that for RBM $X$ in a bounded uniform domain $D\subset \RR^n$,
\begin{equation}\label{eqn:1.1}
\P^x \left( \dim_H X(E)= 2 \dim_H E
\ \hbox{ for all Borel sets } E \subset \RR_+
\right) =1 \quad \hbox{ for every } x \in \overline D.
\end{equation}
Here $X(E)(\omega):=\{X_t(\omega): \, t\in E\}$
denotes the range of $E$ under
RBM $X$ and $\dim_H$ denotes the Hausdorff dimension.
Such a result is called a {\it uniform dimensional result}
because the exceptional set in (\ref{eqn:1.1})
is independent of the Borel time sets
$E\subset \RR_+$. This dependence is important when
one wants to extract information on
$X(E)$ by observing $E$ only
 while $E$ itself is random, for example, when
$E$ is the boundary occupation time set of $X$.

\medskip

We next study the occupation time of $X$ on the boundary $\partial D$,
$$
 S(\omega)=\{ t\geq 0: X_t(\omega) \in \partial D\}.
$$
Corollary \ref{C:3.3} below implies
that
\begin{equation}\label{eqn:1.2}
\dim_H S(\omega ) = 1- \frac{n-\dim_H \partial D }2
\end{equation}
$\P^x$-almost surely for every  $x \in \overline D$.
Note that for any Euclidean domain $D\subset \RR^n$,
$\dim_H \partial D \in [n-1, \, n]$.
This together with (\ref{eqn:1.1}) implies that
the Hausdorff dimension of the boundary trace for RBM in $D$ is
$$\dim_H (X[0,\infty)\cap \partial D) = 2+\dim_H \partial D -n$$
$\P^x$-almost surely for every  $x \in \overline D$.
In particular, for {\it planar} uniform domains
(such as the van Koch snowflake) it follows
that $\P^x$-almost surely
$$ \dim_H (X[0,\infty)\cap \partial D) = \dim_H \partial D
$$
for every  $x \in \overline D$.
This is in contrast to Makarov's celebrated result about the
support of harmonic measure: There is a subset $A$ of $\partial D$
with $\dim_H A = 1$ such that the {\it first} intersection of $X$
with $\partial D$ is almost surely contained in $A$.
In Section \ref{S:5} an example of a (non-uniform) domain
$D$  is given where
the Hausdorff dimension of the boundary trace of RBM in $D$
 {\it differs} from
$\dim_H \partial D$. Thus some assumption about the regularity of
the domain is necessary for our results to hold.

\medskip

The remainder of this paper is organized as follows.
In Section 2, we recall the definitions of RBM and extension domains.
In Section 3, we show that the sample paths of RBM have the same
degree of H\"older continuity as Brownian motion, by
 using the Lyons-Zheng's
forward-backward martingale decomposition of RBM.
We then establish the uniform dimensional results for RBM
on a large class of non-smooth Euclidean domains. For this we derive
the two-sided heat kernel estimates for RBM and use a capacitary
argument.  In the case of dimension
$n=2$, the logarithmic capacity is not good for our approach
so we first subordinate the RBM, then
establish dimensional results for the subordinated stable-like
process and lastly transfer these results back to RBM.
In this procedure, we use the recently obtained
two-sided heat kernel estimates in Chen and Kumagai \cite{CK}
for stable-like processes on $\overline D$.
In Section 4, we apply a recent result in Bogdan, Burdzy and Chen
\cite{BBC} to get the exact
capacitary dimension of the boundary occupation
time of RBM and therefore establish its Hausdorff dimension.
All these results are combined in Section 5 to get the
Hausdorff dimension for the boundary trace of RBM.
Several examples are given to illustrate the main results
of this paper as well as an example that shows that
certain regularity assumptions on the domain $D$ are needed for
the main results of this paper to hold.
Extensions of the main results in this paper to stable-like
jump processes are mentioned in Remarks \ref{R:2.8}, \ref{R:3.4}
and \ref{R:4.2}.

\medskip
We point out here that the approach used
in this paper is quite robust that applies not
only to RBM but also to symmetric reflecting diffusions
as well. These symmetric reflecting diffusions have infinitesimal
generators of divergence form:
$$ {\cal L} =\sum_{i,j=1}^n \frac{\partial}{\partial x_i}
\left( a_{ij}(x) \frac{\partial}{\partial x_j} \right),
$$
where $A(x)=(a_{ij}(x))_{1\leq i, j\leq n}$
is a symmetric $n\times n$ matrix-valued measurable function
that is uniformly elliptic and bounded.
See Chen \cite{C} for relevant informations on symmetric
reflecting diffusions on bounded domains, and
Stroock \cite{St} for two-sided heat kernel estimates
for symmetric diffusions on $\RR^n$.

\section{Extension domains and RBM}
\label{S:2}

Throughout this paper,  $n\geq 2$ is an integer.
For a domain $D\subset \RR^n$, let
$W^{1, 2}(D) :=\{ f\in L^2(D, dx): \, \nabla f \in L^2(D, dx) \}$,
equipped with the Sobolev norm
$\| f \|_{1,2}:= \| f \|_2 + \| \nabla f \|_2$,
where $\| f\|_2:= \left( \int_D f(x)^2 \, dx \right)^{1/2}$.
In this paper we will be concerned with reflecting Brownian motion on domains
$D\subset \RR^n$ whose boundary has zero Lebesgue measure and
have the following $W^{1,2}$-extension property:
there is a bounded linear extension operator
\begin{equation}\label{eqn:extension}
T:  \
W^{1,2}(D) \to
W^{1,2}(\RR^n) \hbox{  such that }
Tf=f  \hbox{ a.e. on } D  \hbox{ for } f\in W^{1,2}(D).
\end{equation}

\medskip

The $(\eps, \delta)$-domains (also called locally uniform domains)
introduced by P. W. Jones \cite{J} have this property.
A domain is $(\eps, \delta)$,
if there are $\eps \in (0, \, \infty)$ and $\delta \in (0, \, \infty]$
such that the conditions in Definition \ref{D:2.1}
for uniform domains hold for $|x-y|<\delta$ with $C=1/\eps$.
Clearly, an $(\eps, \infty)$-domain is just a uniform domain.
For later use, we record the following observation as a lemma.

\medskip

\begin{lemma}\label{L:C1}
If $D\subset \RR^n$ is an $(\eps, \delta)$-domain,
then so is $D\times [0, \, 1]$.
\end{lemma}

\pf
If $(x,t),(y,s)\in D\times [0, \, 1]$
with $|(x,t)-(y,s)|<\text{min}(1, \, \delta )$,
then either $|x-y|\geq |t-s|$ or  $|x-y|< |t-s|$.
In the first case, we can easily lift the
rectifiable curve $\gamma$ from $x$ to $y$ in $D$
in the definition of $(\eps, \delta )$-domain for $D$
to a rectifiable curve from $(x,t)$ to $(y,s)$ that satisfies
the $(\eps, \delta)$-condition for $D\times [0, 1]$,
with possibly different values of $\eps$ and $\delta$.
In the second case, we construct a point $z\in D$ with
$|x-z|=2|t-s|$ such that $\dist (z,\, \partial D))>|t-s|/C$,
and then apply the first case to construct rectifiable
curves joining $(x,t)$
with $(z,(t+s)/2)$ and $(z,(t+s)/2)$ with $(y,s)$.
The piecing together curve satisfies the condition for
$(\eps, \delta )$-domain for $D\times [0, 1]$, with
possibly different values of $\eps $ and $\delta$.
\qed

\medskip

The importance of $(\eps, \delta )$-domains is further
illustrated by the following.
It is shown by Jones \cite{J} that a finitely connected
domain is a $W^{1,2}$-extension domain if and only if
it is an $(\eps, \delta)$-domain. More information on
 domains that have the $W^{1,2}$-extension property
(\ref{eqn:extension}) can be found in Herron and Koskela \cite{HK}.

\medskip

For $f, g \in W^{1,2}(D)$, define
$$\EE(f, g):= \frac12 \int_D \nabla f \cdot \nabla g \, dx,
$$
and $\EE_1 (f, g):= \EE (f, g)+ \int_D f(x) g(x) \, dx$.
We say the Dirichlet form $(\EE, \, W^{1, 2}(D))$
is regular on $\overline D$ if
$W^{1,2}(D)\cap C (\overline D )$ is dense both in
$(W^{1,2}(D),\,  \EE_1^{1/2})$
and in $(C (\overline D), \, \|\cdot \|_\infty )$.
It is known (cf. Theorem 2 on p.14 in Maz'ja \cite{Maz})
that $(\EE, \, W^{1, 2}(D))$
is a regular Dirichlet form on $\overline D$
if $D$ has continuous boundary.
Clearly $(\EE, \, W^{1,2}(\RR^n) )$ is a regular Dirichlet
form on $\RR^n$. So if $D$ is a $W^{1,2}$-extension domain,
$(\EE, \, W^{1, 2}(D))$ is regular on $\overline D$.

\medskip

When $(\EE, \, W^{1, 2}(D))$ is regular on $\overline D$,
there is a strong Markov process $X$ associated with
$(\EE, W^{1, 2}(D))$, having continuous sample paths on $\overline D$;
one can construct a consistent Markovian family of distributions
for the process
starting from every point in $\overline D$ except  possibly for
a subset $N$ of $\partial D$ having zero capacity (see Chen \cite{C}).
When $D$ is smooth, the exceptional set $N$ can be taken to be
the empty set. We will show in Section 3 (see
(\ref{eqn:heat2}) and Lemma \ref{L:2.6} below)
that  the exceptional set $N$ in fact can be taken to
be the empty set for a large class of nonsmooth domains
including $(\eps, \delta)$-domains.
Thus constructed process $X$ is the reflecting
Brownian motion on $D$ in the sense that this definition agrees
with all other standard definitions in smooth domains.
As it is remarked
in Bass, Burdzy and Chen \cite{BBC2}
(see the second paragraph in the proof of
Theorem 5.7 there), such reflecting
Brownian is conformally invariant on planar domains.

\medskip

Let $S(\omega):=\{ t\geq 0: X_t(\omega) \in \partial D\}$
be the occupation time of $X$ on the boundary $\partial D$
and $R(\omega):=X[0, \infty) (\omega)\cap \partial D$
the trace of $X$ on the boundary $\partial D$.

\medskip

Recall that for any increasing function $h$ on $[0, \, 1)$
with $h(0)=0$, one can define a Hausdorff measure ${\cal H}_h$
with respect to the gauge $h$ in the following way (see, e.g.,
p.132 of \cite{AH}).
For $E\subset \RR^n$,
\begin{equation}\label{eqn:haus}
{\cal H}_h (E) = \lim_{\eps \downarrow 0}
 \inf \left\{ \sum_{k=1}^\infty h(r_k): \,
 E \subset \bigcup_{k=1}^\infty B(x_k, r_k) \hbox{ for some }
 x_k\in \RR^n \hbox{ with } \sup_{1\leq k< \infty} r_k \leq \eps
\right\}.
\end{equation}
When $h(r)=r^\alpha$ for some $\alpha>0$, the Hausdorff measure ${\cal H}_h$
will be denoted as ${\cal H}^\alpha$.

\medskip

\begin{defn}\label{def:dset}
A Borel set $\Gamma\subset \RR^n$ is called an $n$-set
if there exists a positive constant $c>0$
such that
$$ m (\Gamma \cap B(x, r))\geq c \,  r^n
\quad \hbox{ for all } x\in \Gamma \hbox{ and } 0<r\leq 1,
$$
where $m$ denotes the Lebesgue measure in $\RR^n$.
\end{defn}

Note that if $\Gamma \subset \RR^n$ is an $n$-set,
then by Proposition 1 in Chapter VIII of \cite{JW},
so is its Euclidean closure $\overline \Gamma$
and $\overline \Gamma \setminus \Gamma$ has
zero Lebesgue measure in $\RR^n$.
It is known that any $(\eps, \delta)$-domain in $\RR^n$ is an
$n$-set (see Example 4 on page 30 of \cite{JW}).

\section{Uniform Dimensional Results for RBM}

In this section, we will extend Kaufman's uniform dimensional
result for planar Brownian motion to reflecting Brownian
motions in $n$-dimensional bounded domains.  This is perhaps
the first time such a uniform dimensional result has been
established for a {\it recurrent} Markov process that does not have
the property of independent increments.

\begin{thm}\label{T:2.1} Suppose that $D\subset \RR^n$ is
a bounded domain whose boundary $\partial D$ has zero Lebesgue
measure such that either
\begin{description}
\item{(1)} it satisfies the
extension property (\ref{eqn:extension}) when $n\geq 3$; or
\item{(2)} $D\subset \RR^2$ is a connected open 2-set such that
 the product domain $D\times [0, 1]$ satisfies
      the extension property (\ref{eqn:extension}).
\end{description}
Then for every $x\in \overline D$,
\begin{equation}\label{eqn:2.1}
\P^x \left( \dim_H X(E)= 2 \dim_H E
\ \hbox{ for all Borel sets } E \subset \RR_+
\right) =1.
\end{equation}
\end{thm}

We will prove the upper bound first.

\begin{thm}\label{T:2.2} Let $D\subset \RR^n$ be a domain
such that $(\EE, W^{1,2}(D))$ is a regular Dirichlet form on $\overline D$.
Then
\begin{description}
\item{(i)} for every $x\in D$,
\begin{equation}\label{eqn:2.2}
\P^x \left( \dim_H X(E)\leq  2 \dim_H E
\ \hbox{ for all Borel sets } E \subset \RR_+
\right) =1;
\end{equation}
\item{(ii)} if $D\subset \RR^n$ satisfies the condition of Theorem
\ref{T:2.1}, then (\ref{eqn:2.2}) holds for every $x\in \overline D$.
\end{description}
\end{thm}

\pf We first show that the sample paths of RBM have the same
degree of H\"older continuity as Brownian motion.
Clearly, since RBM $X_t$ behaves like Brownian motion
when $X_t\in D$, the sample paths of $X$ can not be
smoother than those of Brownian motion.
By Lyons-Zheng's forward and backward martingale decomposition
(see Theorem 5.7.1 of Fukushima, Oshima and Takeda \cite{FOT}),
for any $T>0$,
\begin{equation}\label{eqn:2.3}
 X_t-X_0 = \frac12 W_t - \frac12 \left( W_T
-W_{T-t} \right) \circ r_T
\quad \hbox{for all } 0\leq t\leq T, \ \P^m\hbox{-a.s.},
\end{equation}
where $W$ is a martingale additive functional of $X$
which is an $n$-dimensional Brownian motion,
and $r_T$ is the time reversal operator of
$X$ at time $T$, i.e., $X_t (r_T (\omega ))=X_{T-t}(\omega )$
for each $0\leq t \leq T$.
Since $X$ is symmetric under $\P^m$, $\P^m$ is invariant under
time reversal $r_T$.
For $\alpha \in (0, \, \frac12)$, define
$$ f_\alpha (x) =
\P^x \left( t\to X_t \hbox{ is  $\alpha$-H\"older continuous
   on each finite interval} \right).
$$
It follows from (\ref{eqn:2.3}) that $f_\alpha =1$  q.e. on $ D$.
Here q.e. is the abbreviation for quasi-everywhere; that is
the above property holds holds for every point $x$ in $D$
except a set of zero capacity with respect to the RBM  $X$.
Since $f$ is finely continuous, we have $f_\alpha =1$ q.e. on
$\overline D$. This implies that for q.e. $x\in \overline D$,
$$ \P^x \left( \dim_H X(E)\leq  \alpha^{-1} \,  \dim_H E
\ \hbox{ for all Borel sets } E \subset \RR_+
\right) =1.
$$
Taking a sequence $\alpha_n \uparrow 1/2$, we have
for q.e. $x\in \overline D$,
$$ \P^x \left( \dim_H X(E)\leq  2 \,  \dim_H E
\ \hbox{ for all Borel sets } E \subset \RR_+
\right) =1.
$$
Since $X$ has density function $p(t, x, y)$ for every $x \in D$
(see Fukushima \cite{F}),
the above holds for every $x\in D$.
If $D\subset \RR^n$ satisfies the condition of Theorem \ref{T:2.1},
then by (\ref{eqn:heat2}) and Lemma \ref{L:2.6} below,
$X$ has transition density function $p(t, x, y)$ for
every $x\in \overline D$. It follows that (\ref{eqn:2.2}) holds
for every $x\in \overline D$.
\qed

\medskip

\begin{remark}\label{R:3.3} {\rm
In fact, Theorem \ref{T:2.2}(i) holds for RBM on any Euclidean domain,
with the same proof. See Chen \cite{C} for the construction
of RBM and its forward-backward martingale decomposition
on an arbitrary Euclidean domain.
When $D$ is a domain satisfies the condition of Theorem \ref{T:2.1},
(\ref{eqn:2.2}) can also be established by using the heat kernel
estimates obtained in (\ref{eqn:heat2}) and Lemma \ref{L:2.6}
together with a result initially due to Hawkes and Pruitt
(see Lemma 8.1 of Xiao \cite{Xiao}.
\qed  }
\end{remark}

\medskip

To prove the lower bound in Theorem \ref{T:2.1},
we treat the higher dimensional case ($n\geq 3$)
and two-dimensional case separately.

\medskip

Let $D$ be a domain that satisfies the condition of
Theorem \ref{T:2.1}. Let $N$ be
a Borel set of $\overline D$
having zero capacity with respect to RBM $X$
in $D$
such that  for every $x\in \overline D \setminus N$,
$$\P^x (\hbox{there is some } t>0 \hbox{ such that } X_t \in N
\hbox{ or } X_{t-} \in N)=0.
$$
Such a set $N$  is called a  properly exceptional set of $X$
and alway exists (see \cite{FOT}).
When $n\geq 3$,
by the classical Sobolev inequality on $\RR^n$ and
(\ref{eqn:extension}),
\begin{equation}\label{eqn:2.4}
 \| f \|_{2n/(n-2)} \leq c \, \sqrt{ \EE_1 (f, f)}
\quad \hbox{ for } f \in W^{1,2}(D).
\end{equation}
Thus according to Varopoulos (see Theorem 2.4.2 in Davies \cite{D})
the reflecting Brownian motion $X$ has density function
$p(t, x, y)$ such that
\begin{equation}\label{eqn:heat}
e^{-t} p(t, x, y) \leq c \, t^{-n/2}
\quad \hbox{ for } t>0 \hbox{ and }
x, y \in \overline D \setminus N .
\end{equation}
Recall that it is assumed that
$\partial D$ has zero Lebesgue measure.
Using Davies' method together with an old idea of Nash advanced
by Fabes and Stroock in \cite{FS}, it is now standard
to deduce that (cf. Theorem 2.3 and Theorem 3.4 in \cite{BH}
and Section 3 in \cite{FS}) $p(t, x, y)$ is
jointly H\"older continuous on $\RR_+\times
 D \times D$ and
for any $k>0$ there are positive constants
$c_1, \ c_2 , \ c_3$ and $ c_4 $ such that
\begin{equation}\label{eqn:heat2}
c_1 \,  t^{-n/2} \exp\left( -\frac{|x-y|^2}{c_2 t} \right)
\leq p(t, x, y) \leq
c_3 \,  t^{-n/2} \exp\left( -\frac{|x-y|^2}{c_4 t} \right)
\end{equation}
for $0< t \leq k $ and $x, y \in \overline D$.
So  $p(t, x, y)$ can be extended continuously to
$[0, \infty)\times \overline D \times \overline D$.
This in particular implies that reflecting Brownian motion on
a $W^{1,2}$-extension domain  $D\subset \RR^n$ with $n\geq 3$
can be constructed as a strong Markov process that
starts from {\it every}  point in $\overline D$
(cf. \cite{FOT}).

\medskip

The following covering lemma was first proved for L\'evy processes
in $\RR^n$ by Pruitt and Taylor \cite{PT}. It is extended
to general Markov processes in Liu and Xiao \cite{LX}.

\medskip

\begin{lemma}\label{L:2.3}
Let $\Lambda (a)$ be a collection of cubes of side $a\in (0, \, 1]$
in $\RR^n$  with the property that the number of
these cubes which intersect an arbitrary sphere of radius $a$
in $\RR^n$ is bounded by a constant $K$ that is independent
of $a$ and of the sphere (this happens when the cubes are those
in $a\ZZ^n$ or when the cubes do not
overlap too much). Let $D\subset \RR^n$ satisfy the condition
of Theorem \ref{T:2.1} and let $M(a, t)$ be the number of
those cubes that are hit by RBM $X$ in $D$  before time $t$.
Then there is a constant $c=c(K, t)$ that depends only
on $K$ and $t$ such that
\begin{equation}\label{eqn:2.8}
\E^x \left[ M(a, t) \right] \leq c\, \left[ \inf_{y\in \overline D}
\int_0^t \int_{B(x, a/3)} p(s, y, z) \, dz \, ds \right]^{-1}
\end{equation}
for every $x\in \overline D$.
\end{lemma}

It follows from (\ref{eqn:heat2}) that
\begin{equation}\label{eqn:2.9}
 \E^x \left[ M(a, t) \right] \leq c\, a ^{-2}.
\end{equation}

\begin{lemma}\label{L:2.4}
If $U$ is a ball of radius $a$ in $\RR^n$
and $k>0$,
then there is a constant $c=c(D, k)>0$ such that
$$ \P^x \left( X_s \in U \hbox{ for some } s \in [t, \, k]
\right) \leq c \, \left( \frac{a}{t^{1/2}} \right)^{n-2}
$$
for every $a\leq 1$ and $x\in \overline D$.
\end{lemma}

\pf Let $\mu$ be the 1-equilibrium measure of $X$ for
$U \cap \overline D$; that is,
$$ \E^x [ e^{-\sigma_U} ]=
\int_{\overline D} G_1(x, y) \mu (dy) ,
$$
where $\sigma_U:=\inf \{t>0: \, X_t\in U \}$ and
$G_1(x, y) :=\int_0^\infty e^{-t} p(t, x, y) dt$.
Note that
\begin{eqnarray*}
\P^x \left( X_s \in U \hbox{ for some } s \in [t, \, k] \right)
&\leq & \E^x \left[ \P^{X_t} \left( X_s \in U \hbox{ for some }
             s \in [0, \, k] \right) \right] \\
&\leq & e^k \, \E^x  \left[ \E^{X_t} [ e^{-\sigma_U} ] \right] \\
&\leq & e^k \,  \int_{\overline D} \left(
\int_t^\infty e^{-(s-t)} p(s, x, z) ds \right) \mu (dx) \\
&\leq & c \, e^{2k} \, \int_t^\infty s^{-n/2} ds \, \mu (U) \\
&=& c \, t^{-(n-2)/2} \,  {\rm Cap}_1 (U).
\end{eqnarray*}
It follows from the extension property (\ref{eqn:extension}) that
${\rm Cap}_1 (U) \leq c \, {\rm Cap}_1^{\RR^n} (U)$,
where ${\rm Cap}_1 (U)$ and ${\rm Cap}_1^{\RR^n} (U)$ denote
the 1-capacity of
$U$ for RBM $X$ on $\overline D$
 and Brownian motion in $\RR^n$, respectively.
Hence
$$ {\rm Cap}_1 (U) \leq c \, a^{n-2} ,
$$
which proves the lemma. \qed.

\begin{thm}\label{T:2.5}
Suppose that $n\geq 3$, that $D\subset \RR^n$ is a bounded domain
satisfying the extension property (\ref{eqn:extension}) and that
$\partial D$ has Lebesgue measure.
Then for every $x\in \overline D$,
\begin{equation}\label{eqn:2.11}
\P^x \left( \dim_H X(E)\geq   2 \dim_H E
\ \hbox{ for all Borel sets } E \subset \RR_+
\right) =1.
\end{equation}
\end{thm}

\pf With Lemmas \ref{L:2.3} and \ref{L:2.4} in hand,
by the same argument as those for Lemma 3 and Theorem 1
in Hawkes \cite{H}, we have
$$ \P^x \left( \dim_H X(E\cap [0, \, k])\geq   2 \dim_H
(E \cap [0, \, k])
\ \hbox{ for all Borel sets } E \subset \RR_+
\right) =1
$$
for every $k>0$ and $x\in \overline D$.
Letting $k\to \infty$, we see that (\ref{eqn:2.11}) holds
for every $x\in \overline D $. \qed

\medskip

So far we have proved Theorem \ref{T:2.1}(1). In the remainder of this
section, we will deal with the two-dimensional case.

\begin{lemma}\label{L:2.6}
Suppose $D\subset \RR^2$ is a bounded
domain whose boundary
$\partial D$ has zero Lebesgue measure such that $D\times [0,\,1]$
is a $W^{1,2}$-extension domain in $\RR^3$. Then
reflecting Brownian motion on $D$ exists as a strong Markov
process on $\overline D$ starting from every point in $\overline D$.
Furthermore its transition density function
$p(t, x , y)$ is jointly H\"older continuous on $\RR_+\times
\overline D \times \overline D$ and
has Gaussian estimates (\ref{eqn:heat2}) with $n=2$ there.
\end{lemma}

\pf Under the assumptions of the lemma, $D$ is a $W^{1,2}$-extension
domain in $\RR^2$ so RBM  $X$ on $D$
exists as a strong Markov
process on $\overline D$ starting from every point in
$\overline D$ except on a properly exceptional set $N$ of $X$.
Note that RBM $Y$ in $D\times [0, \, 1]$
can be obtained from $X$ by running an independent
RBM  in the unit interval along the $z$-direction.
By (\ref{eqn:heat}), RBM
$Y$ on $D\times [0,\,1]$ has transition density function
$\widetilde p (t, \wt x, \wt y)$ and
$$ e^{-t} \widetilde p(t, \wt x, \wt y) \leq c \,t^{-3/2}
\quad \hbox{for every } t>0 \hbox{ and }
\wt x, \, \wt y \in (\overline D \setminus N ) \times [0, \, 1].
$$
It follows that $X$ has density function $p(t, x, y)$ and that
$$ e^{-t} p(t, x, y) \leq c \,t^{-1}
\quad \hbox{for every } t>0 \hbox{ and } x , \, y
\in \overline D \setminus N.
$$
Joint H\"older continuity of $p(t, x, y)$ and
its two-sided estimate (\ref{eqn:heat2}) (with $n=2$)
 now follows from a similar argument as that in Section 3 in
Fabes and Stroock \cite{FS} and that for
Theorems 2.3 and 3.4 in Bass and Hsu \cite{BH}.
So the RBM $X$ can be refined to start from every point
in $\overline D$.
 \qed

\begin{remark}\label{R:C2}    {\rm
By Lemma \ref{L:C1} and the remark at the end of Section 2,
any bounded planar $(\eps, \delta )$-domain satisfies the condition
of Lemma \ref{L:2.6}. Combining this with a result of Jones mentioned
at the paragraph following Lemma \ref{L:C1}, we see
that any bounded finitely connected $W^{1,2}$-extension domain in $\RR^2$
satisfies the condition of Lemma \ref{L:2.6}. \qed  }
\end{remark}

\begin{thm}\label{T:2.7} Suppose that $D\subset \RR^2$
is a connected bounded open $2$-set that satisfies
the condition of Lemma \ref{L:2.6}.
Then for every $x\in \overline D$, (\ref{eqn:2.11}) holds.
\end{thm}

\pf Let $X$ be a reflecting Brownian motion in $D$
with sample space $(\Omega, \P^x, x\in \overline D )$.
By Lemma \ref{L:2.6}, the transition density function
has the estimate (\ref{eqn:heat2}) with $n=2$.
For $s \in (0,\, 1)$, let $\xi$ be a $\RR_+$-valued
process on a probability space $(\Omega', \P')$
with independent increments such that
$$ \E \left[ e^{-\lambda \, \xi_t } \right] = e^{-t \lambda^s}
\quad \hbox{for any } \lambda >0
$$
We assume that $X$ and $\xi$ are independent and live on
a common probability space. This can be achieved on
the product space $\Omega \times \Omega$ with
product measures $\{ \P^x \otimes \P', \, x\in \overline D \}$.
Define $Z_t=X_{\xi_t}$; that is, $Z$ is the $s$-subordinate of
$X$.

It follows from Kumagai \cite{K} that $Z$ is a special
case of the stable-like processes on $\overline D$ studied in Chen
and Kumagai \cite{CK}. Hence by Theorem 1.1 of  \cite{CK}
the transition
density function $h(t, x, y)$ of $Y$ has the property that
$$ e^{-t} h(t, x, y) \leq c \, t^{-1/s}
\quad \hbox{ for } t>0 \hbox{ and } x, y \in
\overline D
$$
and
\begin{equation}\label{eqn:2.11b}
 h(t, x, y) \approx
  \min \left\{ t^{-1/s}, \
\frac{t}{|x-y|^{2+2s}} \right\}
\quad \hbox{ for } x, y \in \overline D
\end{equation}
on any finite time intervals.
Here for two function $f, g$, $f \approx g$ means there
are two positive constants $c_1<c_2$ such that
$c_1 g\leq f \leq c_2 g$.
Furthermore, it follows from Theorem 2.5(1) and its proof in
Bogdan, Burdzy and Chen \cite{BBC} that
the capacity induced by $Z$ on $\overline D$
is controlled by the capacity induced by the symmetric
$(2s)$-stable process in $\RR^2$ (cf. also Lemma 3.1 of
Fukushima and Uemura \cite{FU}). Thus by a similar argument
as that for Lemma \ref{L:2.4} above, we have
$$ \P^x \otimes \P'
\left( Z_r \in U \hbox{ for some } r \in [t, \, k]
\right) \leq c \, \left( \frac{a}{t^{1/(2s)}} \right)^{2-2s},
$$
where $U$ is a ball of radius $a$ in $\RR^n$, $k>0$, and
 $c=c(D, k)>0$ is a constant that depends only on $D$, $s$ and $k$.
Using this and applying Lemma \ref{L:2.3} to the process $Z$,
the same reasoning as that in the proof of Theorem \ref{T:2.5}
implies that for every $x\in \overline D$,
\begin{equation}\label{eqn:2.12}
\P^x \otimes \P' \left( \dim_H Z(E)\geq   2s \, \dim_H E
\ \hbox{ for all Borel sets } E \subset \RR_+
\right) =1.
\end{equation}
For a Borel set $E\subset \RR_+$, define
$C_E=\{t\geq 0: \, \xi_t \in E\}$. Clearly
$\dim_H X(E) \geq \dim_H Z(C_E)$,
and the random number
$\dim_H C_E$ depends only on the subordinator $\xi$ and $E$,
which is independent of $X$ and therefore of $\dim_H X(E)$.
Thus by (\ref{eqn:2.12})
and Fubini's theorem, for every $x\in \overline D$,
there is $\Omega_0\subset \Omega$ with $\P^x(\Omega_0)=1$
such that for every $\omega \in \Omega_0$ and for
every Borel measurable set $E\subset \RR_+$,
\begin{equation}\label{eqn:2.13}
\dim_H X(E)(\omega) \geq 2s \, \dim_H C_E (\omega')
\qquad \hbox{for } \P' \hbox{-a.s. } \omega'\in \Omega'.
\end{equation}
Thus for every $\omega \in \Omega_0$ and
every Borel measurable set $E\subset \RR_+$,
$$ \dim_H X(E)(\omega) \geq
2s \, \sup \{ b: \P'( \dim_H C_E > b)>0 \} ,
$$
which by (3.8) of Pruitt \cite{P} is no less than
$2 (s+\dim_H E-1)$.
Hence we have shown that for every $x\in \overline D$,
$$ \P^x \left( \dim_H X(E)\geq   2 (s+ \dim_H E -1)
\ \hbox{ for all Borel sets } E \subset \RR_+
\right) =1.
$$
Letting $s\uparrow 1$ proves the theorem. \qed

\medskip

Theorems \ref{T:2.2}, \ref{T:2.5} amd \ref{T:2.7}
imply Theorem \ref{T:2.1}.

\medskip

\begin{remark}\label{R:2.8}  {\rm
Let $\Gamma \subset \RR^n$ be a closed set such that
there are constants $d\in (0, n]$ and
$c_2>c_1>0$ so that for every $x \in \Gamma$,
$$ {\cal H}^d (B(x, r) \cap \Gamma ) \geq c_1 r^d
 \quad \hbox{for } 0<r\leq 1
\qquad \hbox{and} \qquad {\cal H}^d (B(x, r) \cap \Gamma )
 \leq c_2 r^d
 \quad \hbox{for } r>0
$$
(such a set is called a $d$-set).
Let  $Y$ be an $\alpha$-stable-like
process on $\Gamma$ (see
Chen and Kumagai \cite{CK} for the definition), where
$0<\alpha <2 $ and $\alpha \leq d$.
It is shown in \cite{CK} that $Y$ has a jointly H\"older continuous
transition density function $p(t, x, y)$ on $[0, \infty)\times \Gamma
\times \Gamma$ and that
$$ p(t, x, y) \approx \min \left\{ t^{-d/\alpha}, \
\frac{t}{|x-y|^{d+\alpha}} \right\}
$$
on $[0, k]\times \Gamma \times \Gamma$  for every $k>0$.
Using a similar argument as above,
one can show that
$$
\P^x \left( \dim_H Y(E)\geq  \alpha \dim_H E
\ \hbox{ for all Borel sets } E \subset \RR_+
\right) =1 \quad \hbox{for every } x\in \Gamma.
$$
As it is observed by Xiao in \cite{Xiao},
applying a result initially due to Hawkes and Pruitt (see Lemma 8.1 of
\cite{Xiao}), one gets
$$
\P^x \left( \dim_H Y(E)\leq  \alpha \dim_H E
\ \hbox{ for all Borel sets } E \subset \RR_+
\right) =1 \quad \hbox{for every } x\in \Gamma.
$$
Hence we have for $0<\alpha <2$ and $\alpha \leq d$,
\begin{equation}\label{eqn:2.14}
\P^x \left( \dim_H Y(E) =  \alpha \dim_H E
\ \hbox{ for all Borel sets } E \subset \RR_+
\right) =1 \quad \hbox{for every } x\in \Gamma.
\end{equation}
Taking $E=\RR_+$, one can see clearly that
 the above can not be true when $\alpha>d$.  \qed
}
\end{remark}

\section{Boundary Occupation Time}

\begin{thm}\label{T:3.1}
 Suppose that $D$ is a bounded
 open $n$-set satisfying condition
(\ref{eqn:extension}).
Then for $\gamma \in (0, 1)$,

\begin{description}
\item{(i)} $\dim_H S(\omega) \leq \gamma$
$\P^x$-a.s. for every  $x\in \overline D$ if
$\partial D$ is polar for the symmetric $2(1-\gamma)$-stable
process in $\RR^n$.

\item{(ii)}  $\dim_H S(\omega) \geq \gamma$ $\P^x$-a.s.
for every  $x\in \overline D$
       if $\partial D$
is not polar for the symmetric $2(1-\gamma)$-stable
process in $\RR^n$.
\end{description}
\end{thm}

\pf
For $0<s<1$, let $Z$ be the $s$-subordinator that is independent
of RBM $X$ in $D$.
Then $Y=X_{Z_t}$ is a
recurrent symmetric
Markov process with density functions and its associated
Dirichlet form has the property that its $\EE_1$-norm
is comparable to the Besov norm in $B^{2,2}_s(D)$ (see
Kumagai \cite{K} for a calculation).
Hence by Theorem 2.5 in Bogdan, Burdzy and Chen \cite{BBC},
$\partial D$ is polar for $Y$ if and only if it is polar
for the symmetric $2s$-stable process in $\RR^n$.
Note that according to
Frostman the Hausdorff
dimension is the same as the Riesz capacitary dimension,
while Orey showed that one-dimensional $\alpha$-stable processes
share the same polar sets (see Lemma 2 in Hawkes \cite{H2}).
Therefore we have
for any Borel set $E\subset \RR_+$,
$$\dim_H E = \sup \left\{1-s>0: \ E \hbox{ is not polar
for the } s \hbox{-subordinator } Z \right\} .
$$
So for every $x\in \overline D$, $P^x$-a.s.,
\begin{eqnarray*}
&& \dim_H S(\omega )\\
&=& \sup \left\{ 1-s>0: \ \partial D \hbox{ is not polar
for the $s$-subordinator } Z  \right\}\\
&=& \sup \left\{ 1-s>0: \ \partial D \hbox{ is not polar
for } Y  \right\}\\
&=& \sup \left\{ 1-s>0: \ \partial D \hbox{ is not polar
for the symmetric } (2s) \hbox{-stable process in  } \RR^n \right\} .
\end{eqnarray*}
This proves the theorem. \qed

The above theorem implies that
$\dim_H S(\omega) = \gamma$ almost surely
if and only if $\partial D$ is polar for symmetric $\alpha$-stable
process in $\RR^n$ with $\alpha< 2(1-\gamma)$ and non-polar
for symmetric $\alpha$-stable process with $\alpha> 2(1-\gamma)$.

\begin{remark}\label{R:polar} {\rm
There is
an intimate relationship
between the Hausdorff measure ${\cal H}_h $
(see (\ref{eqn:haus}) for its definition) and the
Riesz capacity ${\rm Cap}_{n-\alpha}$ of order $n-\alpha$,
see Theorems 2.2.7, 5.1.9 and 5.1.13 in \cite{AH}.
Namely,
${\cal H}^{n-\alpha} (A)<\infty$ implies that ${\rm Cap}_{n-\alpha}(A)=0$.
On the other hand if ${\rm Cap}_{n-\alpha} (A)=0$ then
${\cal H}_h (A)=0$
for every $h$ such that
\begin{equation}\label{eqn:Hausdorff}
h \mbox{ is increasing on } [0, \infty) \mbox{ with }
h(0)=0 \quad \mbox{and} \int_0^1 \frac{h(r)} {r^{n+1-\alpha}}
\, dr <\infty.
\end{equation}
In particular,  ${\rm Cap}_{n-\alpha} (A)=0$ implies
 ${\cal H}^\lambda (A)=0$
for any $\lambda>n-\alpha$.
Later we will use the fact that for $n=2$ and $h(t)=1/\log(1/t),$
${\cal H}_h (A)=0$ implies ${\rm cap}\ A=0$,
where ${\rm cap}\ A$ stands  for the logarithmic capacity
see Theorem 5.1.9 in \cite{AH}. \qed
}
\end{remark}

\medskip

Recall that for any Euclidean domain $D\subset \RR^n$,
$\dim_H \partial D \in [n-1, \, n]$.
Combining Remark \ref{R:polar} with Theorem \ref{T:3.1},
 we have the following corollary.

\medskip

\begin{corollary}\label{C:3.3}
Suppose that $D\subset \RR^n$ is a bounded
 open $n$-set satisfying condition
(\ref{eqn:extension}).
Then $\P^x$-a.s.  $\, \dim_H S(\omega ) =
 1- \frac12 \left( n-\dim_H \partial D \right)$
for every $x\in \overline D$.
\end{corollary}

\pf If for some $d\in [n-1, n)$,
${\cal H}^d (\partial D\cap K_m)<\infty $
for an increasing sequence of Borel sets $K_m$
such that  $\cup_{m=1}^\infty K_m\supset \partial D$,
then by Remark \ref{R:polar},
 ${\rm Cap}_d (\partial D)=0$ and so
$\dim_H S(\omega) \leq 1-\frac{n-d}2$ almost surely.
If
${\cal H}^d (\partial D)>0$ for some $d\in [n-1, n)$,
then ${\rm Cap}_d (\partial D)>0$ and hence
 $\dim_H S(\omega ) \geq 1- \frac{n-d}2$ almost surely.
Since
$$ \dim_H \partial D=\inf\{ \alpha>0: \, {\cal H}^\alpha
(\partial D)=0 \} =\sup \{ \alpha>0: \, {\cal H}^\alpha
(\partial D)=\infty \},
$$
the conclusion of the corollary now follows. \qed

\medskip

\begin{remark}\label{R:3.4} {\rm
Let $Y$ be an $\alpha$-stable-like process on an open $n$-set
$D\subset \RR^n$
in the sense of Chen and Kumagai \cite{CK}.
It includes as a special case the
 reflected $\alpha$-stable process on $\overline D$ introduced in
Bogdan, Burdzy and Chen \cite{BBC}.
A similar argument gives the Hausdorff dimension of the
boundary occupation time for $Y$, which asserts that
$\P^x$-almost surely
$$
 \dim_H S^Y(\omega )=\sup \left\{ 1-s>0: \ \partial D \hbox{ is not polar
for symmetric } (\alpha s) \hbox{-stable process in  } \RR^n \right\}
$$
for every $x\in \overline D$.
Hence we have for every $x\in \overline D$, $\P^x$-a.s.
\begin{equation}\label{eqn:stable1}
 \dim_H S^Y (\omega) = \max\left\{ 1 - \frac{n-\dim_H \partial D}{\alpha},
\ 0 \right\} .
\end{equation}
Note that when $\partial D$ has locally finite $d$-dimensional Hausdorff
measure for $d\leq n-\alpha$, it follows from \cite{BBC} that
$\partial D$ is a polar set for $Y$ and therefore $S^Y( \omega)$ is
the empty set almost surely.  \qed
 }
\end{remark}

\section{Boundary Trace of RBM}
\label{S:5}

Combining Theorem \ref{T:2.1} with Corollary \ref{C:3.3}
establishes the following result.

\begin{thm}\label{T:4.1}
Let $D\subset \RR^n$ with $n\geq 2$ be a bounded connected $n$-set
satisfying the condition of
Theorem \ref{T:2.1}.
Then the Hausdorff dimension for the boundary trace
of RBM $X$ in $D$ is
$$ \dim_H \left( X [0, \, \infty ) \cap \partial D \right)
 = 2 +\dim_H \partial D -n \qquad \P^x \hbox{-a.s.}
$$
for every $x\in \overline D$.
\end{thm}

\medskip

\begin{remark}\label{R:4.2} {\rm
Let $D\subset \RR^n$ be an open $n$-set
and $Y$ be an $\alpha$-stable-like
process on $\overline D$  in the sense of \cite{CK}
with $\alpha \leq n$. Then it follows from  Remarks \ref{R:2.8} and
\ref{R:3.4} that the boundary trace $R^Y(\omega):=
Y[0, \infty) \cap \partial D$ of $Y$ has
$$ \dim_H R^Y(\omega)
 =\max \left\{ \alpha +\dim_H \partial D -n, \ 0 \right \}
\qquad \P^x \hbox{-a.s.}
$$
for every $x\in \overline D$.
\qed  }
\end{remark}

\bigskip

\noindent{\bf Example 5.3} Let $D\subset \RR^n$ with $n\geq 2$ be a bounded
Lipschitz domain. Then it is an $n$-set,
satisfies the condition of Theorem \ref{T:2.1}, and
$\dim_H \partial D=n-1$.
So the boundary occupation time set and
boundary trace of RBM in $D$ have Hausdorff dimensions $1/2$ and
$1$, $\P^x$-a.s., respectively,
for every $x\in \overline D$.

\bigskip

\noindent{\bf Example 5.4} Let $D\subset \RR^2$ be
a van Koch snowflake domain.
By Remark \ref{R:C2} it satisfies the condition of
Theorem \ref{T:2.1}(2)
and it is well known that
$\partial D$ has Hausdorff dimension $\frac{\log 4}{\log 3}$.
Hence by Theorem \ref{T:4.1} the boundary occupation time set and
the boundary trace
of RBM in $D$ have Hausdorff dimensions
$\frac12 \frac{\log 4}{\log 3}$ and   $\frac{\log 4}{\log 3} $
$\P^x$-a.s., respectively, for every $x\in \overline D$.

\bigskip

\noindent{\bf Example 5.5} Let $U=D\times (0, 1) \subset \RR^3$,
where $D$ is the van Koch snowflake domain in $\RR^2$. Then
$U$ satisfies the condition of Theorem \ref{T:2.1}(1)
and $\dim_H \partial U =
1+\frac{\log 4}{\log 3}$.
Hence by Theorem \ref{T:4.1} the boundary occupation time set and
the boundary trace
of RBM in $U$ have Hausdorff dimensions
$\frac12 \frac{\log 4}{\log 3}$ and  $\frac{\log 4}{\log 3} $
$\P^x$-a.s., respectively,  for every $x\in \overline U$.
Note that the reflecting Brownian motion in $U$ is $(X, Y)$,
where $X$ is the RBM in $D$ and $Y$ is the RBM in $(0, 1)$.
The above boundary trace result may be a bit surprising if one
compares it with the range and the graph of a 1-dimensional
Brownian motion $B$. It is known that the Hausdorff dimension
for the range $\{ B_t(\omega): \, t\in \RR_+ \}$
of $B$ is 1 a.s., while the Hausdorff dimension for
the graph $\{(t, B_t(\omega )): \, t\in \RR_+\}$ of $B$
is $3/2$ a.s.

\bigskip

\noindent{\bf Example 5.6} Generalizing Example 5.4, let
$D\subset \RR^2$ be a simply connected domain whose boundary is
a regular fractal in the sense of \cite{Mak} (the typical examples
are from complex dynamics, namely components of the Fatou set of a
hyperbolic rational function, and self-similar curves such as the
van Koch snowflake). It is well-known that $\partial D$ has
positive and finite $d$-dimensional Hausdorff measure, where
$d=\dim_H \partial D.$ By Theorem \ref{T:4.1} the boundary
trace of RBM in $D$ has Hausdorff dimension $d$ almost surely.
We will now outline an alternative proof of this result,
based on complex analytic methods.
Let $f:\DD\to D$ be a conformal map of the unit disc
$\DD$ onto $D$.
Using the methods and results from \cite{Mak}, Section 3 (in particular
Proposition 3.2), one can show that there is a set
$A\subset\partial\DD$ with $0<\dim_H A <1$ such that
\begin{equation}\label{eqn:5.1}
\dim_H f(B) = d \qquad \hbox{ for all }
B\subset A  \hbox{ with } \ \dim_H B = \dim_H A.
\end{equation}
In fact, one can show that there is $a\in(0,1)$ such that the set
$$A=\left\{ \xi \in\partial\DD : \,
\lim_{r\to1}\frac{\log |f'(r\xi )|}{|\log(1-r)|}= a \right\}
$$
satisfies $\dim_H f(A)=d$, and that
$f$ is uniformly expanding on $A,$
$$C_{\epsilon}^{-1} |x-y|^{1-a+\epsilon} \leq
 |f(x)-f(y)| \leq C_{\epsilon} |x-y|^{1-a-\epsilon}$$
for every $\epsilon>0$ and all $x,y\in A$.
This easily implies that $A$ satisfies (\ref{eqn:5.1}).
By conformal invariance $f^{-1}(X)$ is
a time change of RBM in $\DD$. Therefore $f^{-1}
\left( X[0, \, \infty) \cap \partial D \right)$
has the same law as the boundary trace of  RBM in $\DD$.
Because RBM in $\DD$ can be constructed from 2-dimensional
Brownian motion by reflecting the part outside $\DD$, and because
Brownian motion intersects $A$ in a set of dimension $\dim_H A$
almost surely by \cite{Per}, we obtain
$\dim_H f^{-1}\left( X[0, \, \infty) \cap \partial D \right)=\dim_H A$
a.s. Now (\ref{eqn:5.1}) implies $\dim_H
\left(  X[0, \, \infty)\cap \partial D \right)
=d$ almost surely.

\medskip

\noindent{\bf Example 5.7}
Our final example shows that some assumption on the regularity of $D$
has to be made in order for Theorem \ref{T:4.1} to hold. We will
construct a simply connected planar domain $D$ such that
$$\partial D = A\cup B,$$
where $A$ has $\sigma$-finite length, $B$ has positive area,
and $f^{-1}(B)$ has zero logarithmic capacity,
where $f: \DD \to D$ is a conformal map. Thus planar Brownian motion
does not visit $f^{-1}(B)$. Using conformal invariance as in
Example 5.6, it follows that $X[0,\infty) \subset D\cup A$ and hence
$\dim_H \left( X[0, \infty) \cap \partial D \right)=1$ a.s., while
$\dim_H \partial D=2$.

\begin{figure}
\centerline{ \scalebox{.6}{\includegraphics{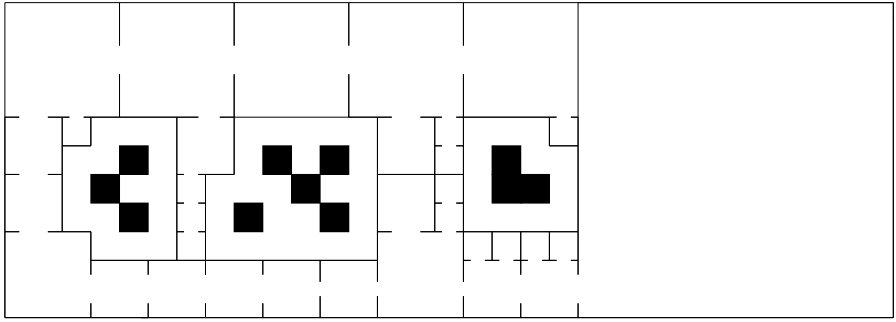}}}
\caption{A finite approximation to B (black squares) and D (white
squares)}
\end{figure}

The idea is to construct $D$
as a countable union of ``smooth'' domains, for instance squares,
joined by narrow ``corridors''. The limit set of the squares is
$B$, and the width of the corridors can be arranged so that
${\rm cap}\ f^{-1}(B)=0$, where ${\rm cap} (A)$ denotes
logarithmic capacity of a set $A \subset \RR^2$.
 To give a rigorous description, it is easier
to begin with the set $B$: Let $B\subset\RR^2$ be a totally disconnected
compact set of positive area in the unit ball $B(0, 1)$.
Consider a decomposition
$$\RR^2 \setminus B = \bigcup_n S_n$$
by closed squares $S_n$ with pairwise disjoint interiors, for instance
the Whitney decomposition of $\RR^2 \setminus B.$ It is easy to find
line segments $I_k$ such that each $I_k$ is contained in some edge of
some $S_n$ and such that
$$D := \left(
\bigcup_n \buildrel \circ \over S_n \cup\  \bigcup_k I_k \right)
\cap B(0, \, 2)
$$
is connected and simply connected, see Figure 1. Then $\partial D = A\cup B$
where $A\subset \left( \cup_n\ \partial S_n \right) \cup
\partial B(0, \, 2)$ has $\sigma$-finite length.
Let $l_k$ denote the length of $I_k$ and $f:\DD\to D$ a conformal
map. Since $\partial D$ is locally connected, $f$ extends continuously
to $\overline \DD$. Given any sequence $\epsilon_k >0,$ it is easy to
choose the $l_k$ such that ${\rm diam} f^{-1}(I_k) < \epsilon_k$
for all $k$ (for instance using Beurling's projection theorem, which gives
${\rm diam} f^{-1}(I_k)\leq C \sqrt{l_k}$).
Denote $J_k$ the arc on $\partial\DD$ that is separated
from 0 by $f^{-1}(I_k)$, so that ${\rm diam} J_k \leq 2 \epsilon_k.$
If $x\in f^{-1}(B),$ then $f[0,x)$ passes through infinitely many
$I_k.$ Therefore $f^{-1}(B)\subset \cup_{n\geq k} J_n$ for all $k.$
By choosing $\epsilon_k$ such that
$\sum_{n=1}^{\infty} 1/\log(1/\epsilon_n) < \infty$ we get
${\cal H}_h (f^{-1}(B))=0$ with $h(t)=1/\log(1/t).$
Now ${\rm cap}\ f^{-1}(B)=0$ follows from Remark \ref{R:polar}.

\bigskip
\bigskip

\noindent{\bf Acknowledgements}. We thank Y. Xiao for very
helpful comments on a preliminary version of this paper.

\vskip 0.5truein

\begin{singlespace}

\end{singlespace}

\vskip 0.3truein

Itai Benjamini

Weizmann Institute, Rehovot 76100, Israel

Email: {\tt itai@wisdom.weizmann.ac.il}

\bigskip

Zhen-Qing Chen \  and \ Steffen Rohde

Department of Mathematics, University of Washington, Seattle,
WA 98195, USA

E-mail: {\tt zchen@math.washington.edu} \  and \
 {\tt rohde@math.washington.edu}

\end{doublespace}

\end{document}